\documentclass[leqno,draft]{article}

\begin{document}

\title{Some aspects of calculus on non-smooth sets}

\author{Stephen Semmes \\
        Rice University}

\date{}

\maketitle

        Let $E$ be a closed set in ${\bf R}^n$, and suppose that there
is a $k \ge 1$ such that every $x, y \in E$ can be connected by a
rectifiable path in $E$ with length $\le k \, |x - y|$.  This
condition is satisfied by chord-arc curves, Lipschitz manifolds of any
dimension, and fractals like Sierpinski gaskets and carpets.  Note
that length-minimizing paths in $E$ are chord-arc curves with constant
$k$.

        A basic feature of this condition is that one can integrate
local Lipschitz conditions on $E$ to get global conditions.  For
instance, if $f : E \to {\bf R}$ is locally Lipschitz of order $1$
with constant $C \ge 0$, then $f$ is globally Lipschitz on $E$ with
constant $k \, C$.

        Let $A(x)$ be a continuous function on $E$ with values in
linear mappings from ${\bf R}^n$ to ${\bf R}$.  Equivalently, one can
use the standard inner product on ${\bf R}^n$ to represent $A(x)$ by a
continuous mapping from $E$ into ${\bf R}^n$.  Also let $f$ be a
locally Lipschitz real-valued function on $E$.

        Suppose that $A$ includes the directional derivatives of $f$
almost everywhere on any rectifiable curve in $E$, in the sense that
the derivative of $f(p(t))$ is equal to $A(p(t))$ applied to the
derivative of $p(t)$ for almost every $t$ when $p(t)$ is a locally
Lipschitz function on an interval $I$ in the real line with values in $E$.
In particular, this holds when $f$ is continuously differentiable on
$E$ in the sense of the Whitney extension theorem with differential at
$x \in E$ given by $A(x)$.  In this case,
\begin{equation}
        f(p(b)) - f(p(a)) = \int_a^b A(p(t))(p'(t)) \, dt
\end{equation}
for every $a, b \in I$.  If $A(x)$ is constant on $E$, then it follows
that $f$ is equal to the restriction of an affine function on ${\bf
R}^n$ to $E$.

        If $A$ is not constant, then we still have that
\begin{eqnarray}
\lefteqn{|f(y) - f(x) - A(x)(y - x)|} \\
 & & \le k \, |x - y| \, \sup \{|A(w) - A(x)| :
                    w \in E, |x - w| \le k \, |x - y|\} \nonumber
\end{eqnarray}
for every $x, y \in E$.  This follows from the previous formula
applied to a curve in $E$ that connects $x$ to $y$ and has length $\le
k \, |x - y|$, by approximating $A(p(t))$ in the integral by $A(x)$
and estimating the remainder by the oscillation of $A$ on the path
times the length of the path.  If we did not know already that $f$ is
continuously differentiable on $E$ in Whitney's sense, then it follows
from this and the continuity of $A$.  Even if we did know that $f$ is
continuously-differentiable on $E$, this estimate provides more
precise information about the behavior of $f$ in terms of the
continuity properties of $A$.  For example, H\"older continuity of
order $\alpha$ of $A$ implies $C^{1, \alpha}$-type smoothness of $f$.

        If $E$ is a Lipschitz submanifold of ${\bf R}^n$, then $E$ has
tangent planes almost everywhere and a locally Lipschitz function $f$
on $E$ has linear differentials on almost all of these tangent planes.
If there is a continuous function $A(x)$ on $E$ with values in linear
functionals on ${\bf R}^n$ such that these differentials of $f$ are
equal to the restriction of $A(x)$ to the tangent plane to $E$ at
almost everywhere $x \in E$, then similar reasoning shows that $f$ is
continuously differentiable on $E$ in Whitney's sense.

        If $E$ happens to be a $C^1$ submanifold of ${\bf R}^n$, then
$f$ is a $C^1$ function on $E$ in the usual sense.  Otherwise,
restrictions of smooth functions on ${\bf R}^n$ to $E$ may not be
smooth in a more intrinsic way.  It may be that $E$ is a non-smooth
embedding of a smooth curve or surface, and that smooth functions on
${\bf R}^n$ do not correspond to smooth functions on the parameter
domain.  It may be that $E$ is something like a Lipschitz manifold for
which there is not a compatible smooth structure.  Nonetheless, we get
some regularity to extrinsic smoothness under the geometric condition
under consideration.

        If $E$ is a chord-arc curve, then one can integrate a
continuous family $A(x)$ of differentials on $E$ to get a function
$f$.  This is analogous to integrating a continuous function on an
interval to get a $C^1$ function, but it is not quite the same, since
the smoothness is now defined extrinsically.

        Let us emphasize that this type of geometric condition does
not imply anything like Poincar\'e or Sobolev inequalities, aside from
the special case of chord-arc curves.  For instance, one can have
bubbles with small boundaries where a curve can easily pass.  Thus one
works with quantities based on the supremum norm instead of integrals.

        Just as classical Fourier analysis on the real line deals with
the ordinary differential operator $d/dx$, one can look at analysis on
Lipschitz graphs in the complex plane as dealing with certain
perturbations of $d/dx$, as in \cite{c-m-1, c-m-2}.  The complex
derivative of a holomorphic function on a neighborhood of the graph at
a point on the graph is a multiple of the ordinary derivative tangent
to the graph, which corresponds to a perturbation of $d/dx$ when the
graph is parameterized by projecting to the real line.  This also
works for a complex-valued continuously-differentiable function on the
graph in Whitney's sense for which the differentials are
complex-linear.

        Additional conditions for the differentials of
continuously-differentiable functions in Whitney's sense like
complex-linearity can be very interesting.  In particular, they may
imply uniqueness of the differentials.  As another version of this, a
real-linear mapping from ${\bf R}^n$ into a Clifford algebra is
uniquely determined by its restriction to any hyperplane when it is
left or right Clifford holomorphic.  Thus it is natural to look at
Clifford-valued continuously-differentiable functions on hypersurfaces
in ${\bf R}^n$ for which the differentials are left or right Clifford
holomorphic.  Note that the differentials of a
continuously-differentiable function in Whitney's sense are uniquely
determined by the function at any point where the set is not
approximately contained in a hyperplane.

        Normally, the differentials of a continuously-differentiable
function in Whitney's sense are not uniquely determined by the
function on the set, and they may not be very stable even if they are
uniquely determined.  The set may be contained in a smooth submanifold
of lower dimension, so that the behavior of the differentials in the
normal directions is not important.  If there is suitable uniqueness
and stability, then continuity properties of the differentials can be
derived from the corresponding approximation properties of the
function by affine functions on the set.

        Uniqueness of the differential, perhaps through an additional
condition like complex-linearity, has the effect of allowing a
derivative in any direction to be interpreted as being tangent to the
set.  In this way, extrinsic behavior becomes more intrinsic.

\end{document}